\documentclass[final,1p,times]{elsarticle}
\usepackage{amssymb}
\usepackage{amssymb}
\usepackage{amsmath}
\usepackage{amsthm}
\usepackage[utf8]{inputenc}
\usepackage{marginnote}
\usepackage{amsfonts}
\usepackage{textcomp}
\usepackage{color,soul}
\usepackage{placeins}
\usepackage{tikz}
\usepackage{relsize}
\usepackage{makecell}
\usepackage{algorithm}
\usepackage{algpseudocode}
\usepackage{amssymb}
\usepackage{graphicx}
\usepackage{subcaption}
\usepackage{hyperref}
\usepackage[none]{hyphenat}
\usepackage{lipsum}
\usepackage{tikz}

\newdefinition{defn}[thm]{Definition}
\allowdisplaybreaks
\journal{}

\author[sha]{Shankhadeep Mondal}
\ead{shankhadeep.mondal@ucf.edu}
\author[ram]{R. N. Mohapatra}
\ead{ram.mohapatra@ucf.edu}

\address{Department of Mathematics, University of Central Florida}
		\cortext[shan]{Corresponding Author: shankhadeep.mondal@ucf.edu}

\begin{document}
\begin{frontmatter}

\title{When Mathematics Meets Painting: Fibonacci Geometry, Cubism and Visual Abstraction}

\begin{abstract}
This paper explores the Fibonacci sequence and the Golden Ratio as organizing principles for visual composition and abstraction in painting. The author shows how recursive proportional systems, long associated with natural growth and aesthetic harmony, inform artistic structure and visual balance. The discussion traces Fibonacci-based geometry from Renaissance art to modern and contemporary practices, with particular attention to Cubism, where fragmentation and multiple viewpoints echo principles of recursion and geometric division. Through selected artistic examples and mathematical insight, the paper demonstrates that Fibonacci geometry functions not merely as a symbolic reference but as a generative framework shaping visual abstraction and artistic expression.
\end{abstract}

\begin{keyword}
Fibonacci sequence \sep Golden Ratio \sep Mathematical art \sep Visual composition,  Cubism 
\end{keyword}

\end{frontmatter}

\section{Introduction}

The dialogue between mathematics and visual art has shaped artistic theory and practice for centuries, from classical studies of proportion to modern explorations of abstraction and form. Mathematical concepts such as symmetry, geometry, and numerical ratio have long served as organizing principles for visual space and aesthetic balance. Among these, the Fibonacci sequence holds a distinctive position due to its recursive structure and its close connection to the Golden Ratio, a proportion historically associated with harmony, growth, and visual coherence \cite{livio2002golden, ghyka1977geometry}. The Fibonacci sequence $\{x_n\}_{n \geq 0}$ is defined by the recurrence relation
\[
x_{n+2} = x_{n+1} + x_n, \qquad x_0 = 0,\; x_1 = 1,
\]
and generates ratios of successive terms that converge to the Golden Ratio
\[
\varphi = \frac{1+\sqrt{5}}{2}.
\]
This convergence links discrete recursion with continuous proportionality and appears widely in natural systems, including phyllotaxis, branching patterns, and biological growth \cite{huntley1972divine, doczi1994power}. Such visual regularities have inspired artists to employ Fibonacci-based grids, spirals, and scaling principles as compositional frameworks reflecting both natural order and mathematical structure.

During the Renaissance, artists and theorists emphasized geometry, proportion, and perspective as foundations of visual representation. Although Fibonacci numbers were rarely cited explicitly, proportional systems closely related to the Golden Ratio informed studies of anatomy, architecture, and pictorial balance \cite{field1997invention}. Mathematics thus functioned as an underlying structural guide rather than a decorative element.

In the twentieth century, Cubism introduced a decisive transformation in the mathematical organization of visual form. By rejecting single-point perspective, Cubist artists fragmented objects into geometric facets and reassembled them through multiple viewpoints. This process may be interpreted mathematically as geometric decomposition and recombination, governed by proportional division, recursive subdivision, and modular organization \cite{henderson1983fourth}. While Cubist painters did not explicitly invoke Fibonacci theory, the internal logic of Cubism closely aligns with Fibonacci principles such as repetition with variation, hierarchical scaling, and balanced part–whole relationships.

The aim of this paper is to examine how the Fibonacci sequence functions as a unifying mathematical framework across representational, Cubist, and abstract painting. Rather than treating Fibonacci patterns as decorative motifs, we emphasize their structural role in organizing visual space and guiding abstraction. Through mathematical analysis, historical context, and selected artistic case studies, we demonstrate how Fibonacci geometry continues to inform contemporary visual abstraction and mathematically grounded painting.

\section{Historical Background}

The relationship between mathematics and art has deep historical roots, extending back to classical antiquity and reaching a mature synthesis during the Renaissance. Artists and scholars of this period increasingly recognized that numerical proportion, geometric harmony, and spatial order could serve as foundations for visual beauty. Among these mathematical tools, the Fibonacci sequence was introduced to the Western world by Leonardo of Pisa (Fibonacci) in the 13th century and emerged as a key structure linking natural growth patterns to aesthetic principles. 

Renaissance artists such as Leonardo da Vinci explored proportional systems derived from mathematics, including the golden ratio, which is intimately connected to the Fibonacci sequence. Although Leonardo did not explicitly reference Fibonacci numbers, his studies on anatomical proportion, architectural harmony, and geometric form reveal a philosophical commitment to numerical balance and modular growth that resonates with Fibonacci’s recursive structure. This emphasis on hidden order and proportional logic would later reappear in modern art, where classical perspective was deliberately dismantled in favor of structural analysis.

In the twentieth century, Cubism offered a radical reinterpretation of mathematical structure in art by fragmenting objects into geometric planes and reassembling them according to internal compositional logic rather than optical realism. While Cubist artists such as Pablo Picasso and Georges Braque did not explicitly employ the Fibonacci sequence, their emphasis on modular decomposition, repetition, and proportional balance parallels the recursive and hierarchical principles underlying Fibonacci growth. These ideas were further abstracted by artists such as Wassily Kandinsky, who explicitly associated form, geometry, and rhythm with mathematical and spiritual order. Kandinsky’s abstract compositions emphasize dynamic balance, geometric relationships, and structural harmony, aligning closely with Fibonacci-inspired notions of progression and self-similarity. In parallel, artists such as Salvador Dalí directly incorporated the golden ratio and Fibonacci-based proportions into works like \emph{The Sacrament of the Last Supper}, demonstrating an explicit synthesis of mathematical ratio and artistic intention.

This historical trajectory from Renaissance proportion, through Cubist structural analysis, to Kandinsky’s abstraction and deliberate Fibonacci-based compositions provides the foundation for contemporary explorations of the Fibonacci sequence in painting. It illustrates that the dialogue between mathematics and visual art extends beyond explicit numerical usage, encompassing deeper structural principles that continue to evolve across artistic movements and aesthetic frameworks.

\section{Mathematical Foundations of the Fibonacci Sequence}
The Fibonacci sequence is defined as a series of numbers where each number is the sum of the two preceding ones. It begins with 0 and 1, and continues indefinitely: 0, 1, 1, 2, 3, 5, 8, 13, and so on. The ratio of successive Fibonacci numbers approximates the golden ratio (approximately 1.618), a number that has been associated with aesthetic harmony and beauty.\\

The Fibonacci sequence, defined by the recurrence relation
\[
F_{n+2} = F_{n+1} + F_{n}, \qquad F_0 = 0, \; F_1 = 1,
\]
is one of the most elegant and pervasive structures in mathematics and nature. Its recursive structure embodies the concept of self-generation and growth, reflecting balance between simplicity and complexity — an idea that resonates profoundly with visual composition in mathematical painting.

\subsection{Closed Form and the Golden Ratio}

The Fibonacci sequence admits a closed-form expression, known as \emph{Binet’s formula}:
\[
F_n = \frac{\varphi^n - \psi^n}{\varphi - \psi},
\]
where $\varphi = \frac{1+\sqrt{5}}{2}$ is the \textit{Golden Ratio} and $\psi = \frac{1-\sqrt{5}}{2}$ is its conjugate. The ratio of successive Fibonacci numbers satisfies
\[
\lim_{n \to \infty} \frac{F_{n+1}}{F_n} = \varphi.
\]
This convergence encapsulates the aesthetic ideal of proportion and harmony. In visual art and geometry, $\varphi$ governs spatial balance — the placement of elements such that neither dominates the other, creating a dynamic equilibrium that appeals naturally to human perception.

\subsection{Geometric and Analytical Interpretations}

Each Fibonacci number can be geometrically represented as a tiling or partition of space. For example, the \emph{Fibonacci spiral} is obtained by connecting quarter circles inside squares with side lengths equal to Fibonacci numbers. The resulting spiral approximates a logarithmic spiral, expressed in polar coordinates as
\[
r = a e^{b\theta}, \quad \text{where}\quad b = \frac{\ln \varphi}{\pi/2}.
\]
This curve symbolizes self-similar growth and appears in diverse natural phenomena from sunflower seed arrangements to the shells of nautilus and the branching of trees. In mathematical art, this spiral captures the idea of infinite expansion from finite recursion, bridging algebraic recurrence and geometric beauty.

\subsection{Applications and Mathematical Connections}

The Fibonacci sequence connects deeply with multiple branches of mathematics:

\begin{itemize}
\item \textbf{Number theory:} Fibonacci numbers satisfy elegant divisibility and modularity properties. For instance, $\gcd(F_m,F_n)=F_{\gcd(m,n)}$, and the sequence modulo $k$ is always periodic (Pisano period).
\item \textbf{Combinatorics:} $F_n$ counts the number of ways to tile a $1 \times (n-1)$ board using $1\times1$ and $1\times2$ tiles. It also enumerates binary strings of length $n-1$ with no consecutive ones.
\item \textbf{Linear algebra:} The recurrence relation can be expressed as a matrix power:
\[
\begin{pmatrix}
F_{n+1}\\ F_n
\end{pmatrix}
= 
\begin{pmatrix}
1 & 1\\
1 & 0
\end{pmatrix}^n
\begin{pmatrix}
1\\ 0
\end{pmatrix}.
\]
This form enables spectral analysis; the eigenvalues of the matrix are $\varphi$ and $\psi$, reaffirming the central role of the Golden Ratio.

\item \textbf{Fractals and dynamical systems:} 
The Fibonacci recurrence arises naturally in the study of 
substitution dynamical systems. Consider the Fibonacci substitution
\[
\sigma:\quad a \mapsto ab,\qquad b \mapsto a.
\]
The associated substitution matrix is
\[
M_\sigma=\begin{pmatrix}
1 & 1\\
1 & 0
\end{pmatrix},
\]
whose eigenvalues are $\varphi$ and $\psi$.  
Iterating $\sigma$ generates infinite words whose letter frequencies converge to the normalized Perron--Frobenius eigenvector of $M_\sigma$, giving the limiting ratio $\varphi$ of $a$'s to $b$'s. Geometrically, coding these words by intervals produces the
\emph{Fibonacci quasicrystal}, a one-dimensional self-similar tiling
with inflation factor $\varphi$ satisfying
\[
\Lambda = \varphi \Lambda \,\cup\, (\varphi \Lambda + 1),
\]
a functional equation expressing exact self-similarity.  
The spectral properties of the tiling dynamical system $(X_\sigma, T)$ are governed by $M_\sigma$, linking Fibonacci growth to the pure point diffraction spectrum characteristic of quasi-periodic structures such as Penrose tilings.

\end{itemize}

\subsection{Artistic Illustrations: From Leonardo to Turner}

Mathematical harmony, as embodied by the Fibonacci sequence and the Golden Ratio, has long served as a guiding principle in visual art. The interplay between proportion, symmetry, and visual balance can be vividly observed in the compositions of both Renaissance and Romantic masters. This section highlights how mathematical ideas are manifested in two celebrated examples: Leonardo da Vinci’s \emph{Mona Lisa} and J. M. W. Turner’s landscape paintings.

\subsubsection*{Leonardo da Vinci’s \emph{Mona Lisa} and the Golden Ratio}

Leonardo’s \emph{Mona Lisa} (c.~1503–1517), housed in the Louvre Museum, is often cited as a paradigm of mathematical beauty in art. Numerous geometric analyses reveal that the spatial organization of the portrait approximates the Golden Ratio $\varphi \approx 1.618$, which is directly related to the Fibonacci sequence by
\[
\lim_{n \to \infty} \frac{F_{n+1}}{F_n} = \varphi.
\]
When a golden rectangle or Fibonacci spiral is superimposed on the painting, key features such as the head, eyes, and hands are align closely with its curve. This suggests that the composition naturally echoes recursive proportionality, where each visual segment grows in harmonious ratio to the preceding one.

\begin{figure}[h!]
\centering
\includegraphics[width=0.45\textwidth]{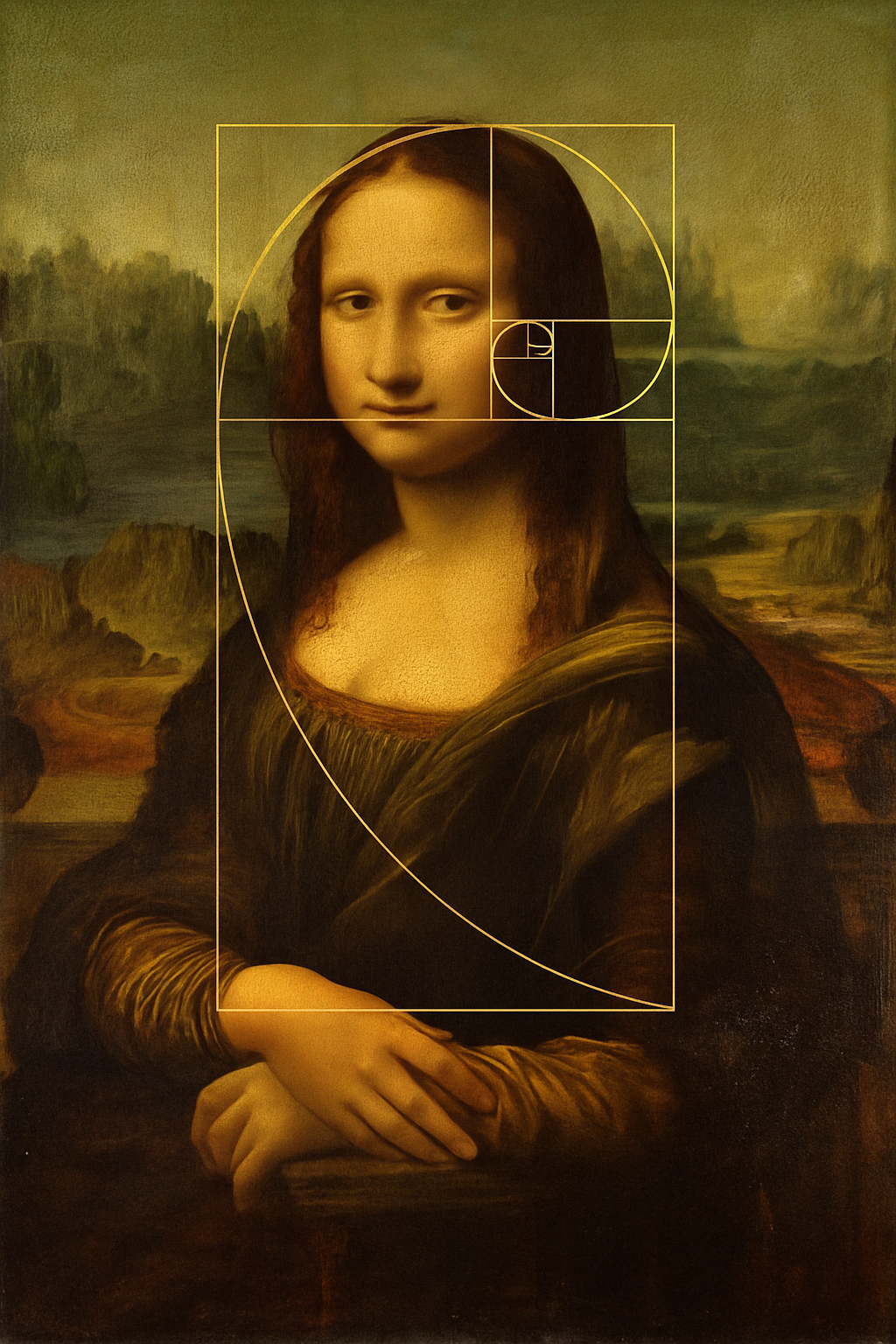}
\caption{Figure 4.2: Leonardo da Vinci’s \emph{Mona Lisa}, often analyzed as a geometric embodiment of the Golden Ratio. The overlay of Fibonacci spirals and golden rectangles demonstrates near ideal proportional balance between the figure’s facial region and torso.}
\end{figure}

Although there is no definitive historical evidence that Leonardo consciously applied the Golden Ratio formula, his deep study of geometry, anatomy, and proportion recorded in the \emph{Codex Atlanticus} suggests a strong awareness of mathematical harmony. The Fibonacci-based spiral division in \emph{Mona Lisa} reflects an organic equilibrium: a visual metaphor for self-similar growth and the recursive perfection of natural forms.

\subsubsection*{J.~M.~W. Turner and the Geometry of Light}

Joseph Mallord William Turner (1775–1851), one of Britain’s most visionary Romantic painters, integrated mathematics, geometry, and scientific observation into his practice. Appointed \emph{Professor of Perspective} at the Royal Academy, Turner explored geometric projection, vanishing points, and atmospheric optics with mathematical precision. His sketchbooks contain constructions of circles, spirals, and linear diagrams used to guide the composition of his luminous seascapes and skies.

Two particular works exemplify this union of art and mathematics:

\begin{enumerate}
\item \textbf{\emph{Norham Castle, Sunrise}} (c.~1845): Turner structured the composition through linear perspective and horizon geometry. The placement of the castle and reflection along horizontal $\varphi$-divisions of the canvas creates a visual equilibrium analogous to the Golden Ratio. The gradation of color intensity across the sky can be interpreted as a continuous exponential decay which is  mathematically expressed as $I(x) = I_0 e^{-kx},$ representing light diffusion.

\begin{figure}[h!]
\centering
\includegraphics[width=0.70\textwidth]{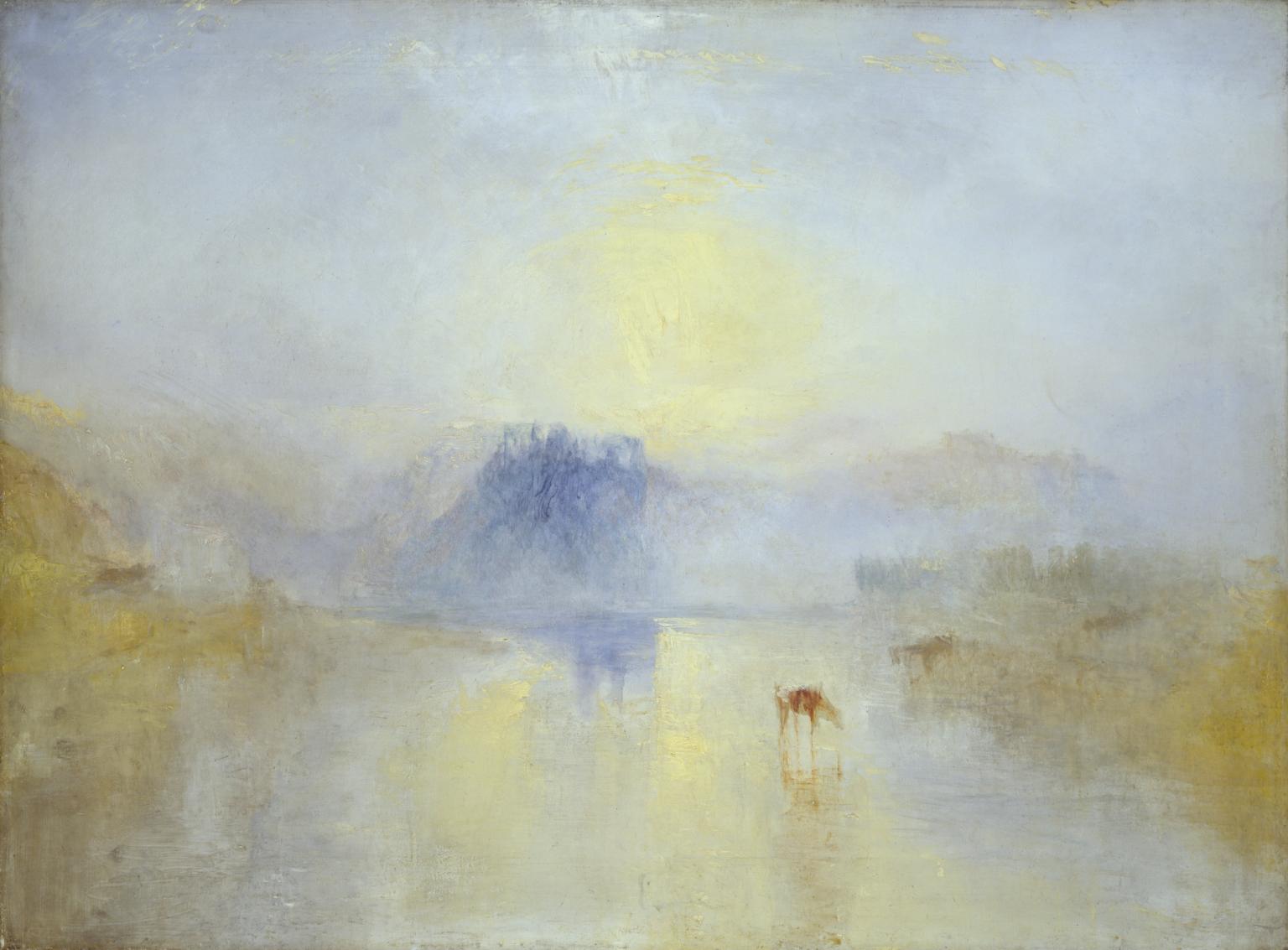}
\caption{Figure 4.3: J.~M.~W. Turner’s \emph{Norham Castle, Sunrise}. The geometric layout of the horizon and castle placement corresponds approximately to Golden Ratio divisions, while the gradient of luminance follows an exponential decay pattern, illustrating the mathematical structure underlying atmospheric perspective.}
\end{figure}

\item \textbf{\emph{Frosty Morning}} (1813): This work demonstrates precise control of vanishing lines, leading the viewer’s eye toward a focal point in the distance. The spatial contraction obeys the principle of geometric projection, where parallel lines meet at infinity—an elegant manifestation of projective geometry in art.\\
Turner’s \emph{Frosty Morning} reveals a tightly organized geometric structure shaped by Golden Ratio divisions, logarithmic spirals, and classical compositional grids. The vertical and horizontal partitions at $x = W/\varphi$ and $y = H/\varphi$ align closely with major narrative elements particularly the horses, cart, and luminous horizon—thus dividing the landscape into harmonically balanced regions. A sequence of nested Golden rectangles of widths $W,\, W/\varphi,\, W/\varphi^2,\ldots$ encloses progressively finer visual features, forming a recursive $\varphi$-scaled hierarchy that echoes the painting’s rhythmic spatial flow. Superimposing an approximate Fibonacci (logarithmic) spiral, generated by quarter-circle arcs satisfying $r(\theta)=ae^{b\theta}$ with $b=\ln(\varphi)/(\pi/2)$, traces a natural viewing trajectory from the foreground through the central horse and cart to the distant atmospheric light. The canvas diagonals intersect near the horse, cart wheel, and right-hand worker, reinforcing the compositional center, while the rule-of-thirds grid coincides with the placements of the tree line, family group, and laborer. Together, these Golden Ratio lines, spiral curves, diagonals, and third divisions concentrate visual attention around key focal points, revealing a compact yet powerful mathematical framework underlying Turner’s atmospheric landscape.

\begin{figure}[h!]
\centering
\includegraphics[width=0.70\textwidth]{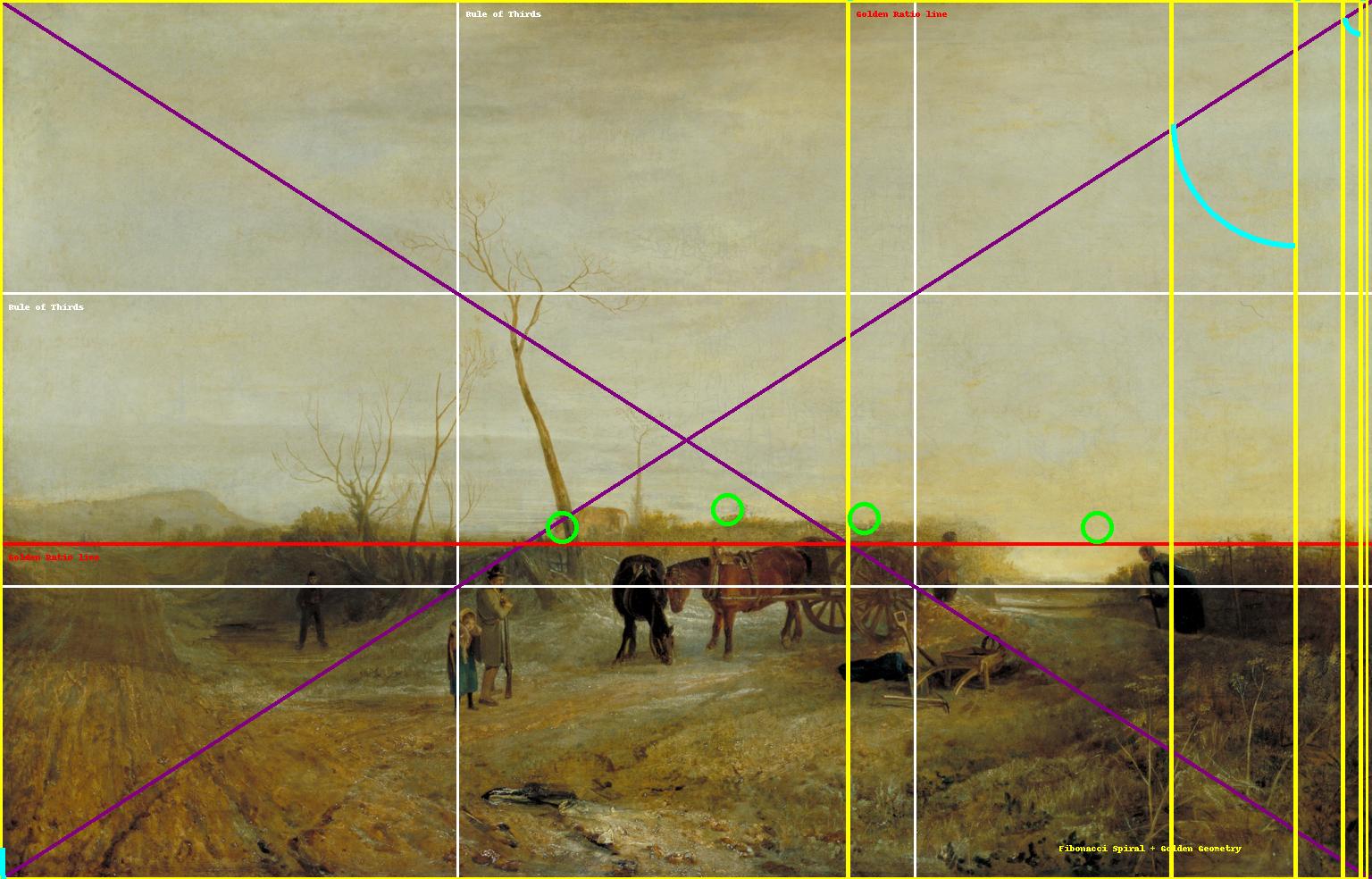}
\caption{Figure 4.4: J.~M.~W. Turner’s \emph{Frosty Morning}. 
The arrangement of the horizon, road, and architectural forms exhibits proportional divisions close to the Golden Ratio, establishing a balanced visual geometry. Moreover, the distribution of luminance—moving from the bright, mist-filled foreground to the muted, atmospheric distance—follows an approximate exponential decay profile. Together, these structural features reveal the mathematical organization underpinning Turner’s treatment of atmospheric depth and early-morning illumination.}
\end{figure}

\end{enumerate}
Through these works, Turner bridges science and art—demonstrating that the play of light, space, and form can be grounded in quantitative principles. His use of perspective and luminance aligns with mathematical models of reflection and decay, turning empirical observation into visual poetry.

\bigskip

Thus, from Leonardo’s proportional geometry to Turner’s spectral perspective, mathematics provides both the invisible architecture and the philosophical foundation of art. In the language of Fibonacci growth and geometric harmony, these masterpieces remind us that beauty often arises where mathematical order and creative intuition intersect.

\subsection{Kandinsky and Mathematical Abstraction}

Wassily Kandinsky (1866--1944), often regarded as the pioneer of abstract art, developed a visual 
language grounded in mathematical structure. In his treatise \emph{Point and Line to Plane} (1926), 
Kandinsky formalized painting as a geometric system: a point is a ``zero-dimensional nucleus,'' a line 
is the trajectory of a point under motion, and a plane is the field in which these elements interact.  
This framework mirrors fundamental concepts in geometry and dynamical systems, where motion generates 
higher-dimensional objects.  
Kandinsky's compositions frequently employ symmetries, rotations, and affine transformations, creating 
visual analogues of linear operators acting on geometric primitives.  
His recursive use of circles, triangles, and arcs produces quasi-fractal layering, reinforcing the 
structural logic seen in Fibonacci spirals and self-similar growth.  
Thus, Kandinsky’s work forms a bridge between abstract painting and mathematical thinking, revealing an 
aesthetic grounded in geometry, proportion, and the algebra of visual forms.

\subsubsection{Mathematical Analysis of Kandinsky’s \emph{Composition VII} and \emph{Composition VIII}}

Kandinsky’s \emph{Composition VII} (1913) and \emph{Composition VIII} (1923) represent two distinct yet
mathematically related phases of abstraction.
While \emph{Composition VII} exhibits continuous, fluid, and highly nonlinear motion fields,
\emph{Composition VIII} is governed by rigid geometric primitives such as circles, triangles, and lines,
arranged according to an underlying operator system.
Both works may be interpreted in the Hilbert space
\( H = L^{2}(\Omega; \mathbb{R}^{3}) \), where each visual element is generated by the action of geometric
operators on a finite set of basis shapes.

\begin{figure}[h!]
\centering

\begin{subfigure}[t]{0.45\textwidth}
    \centering
    \includegraphics[width=\textwidth]{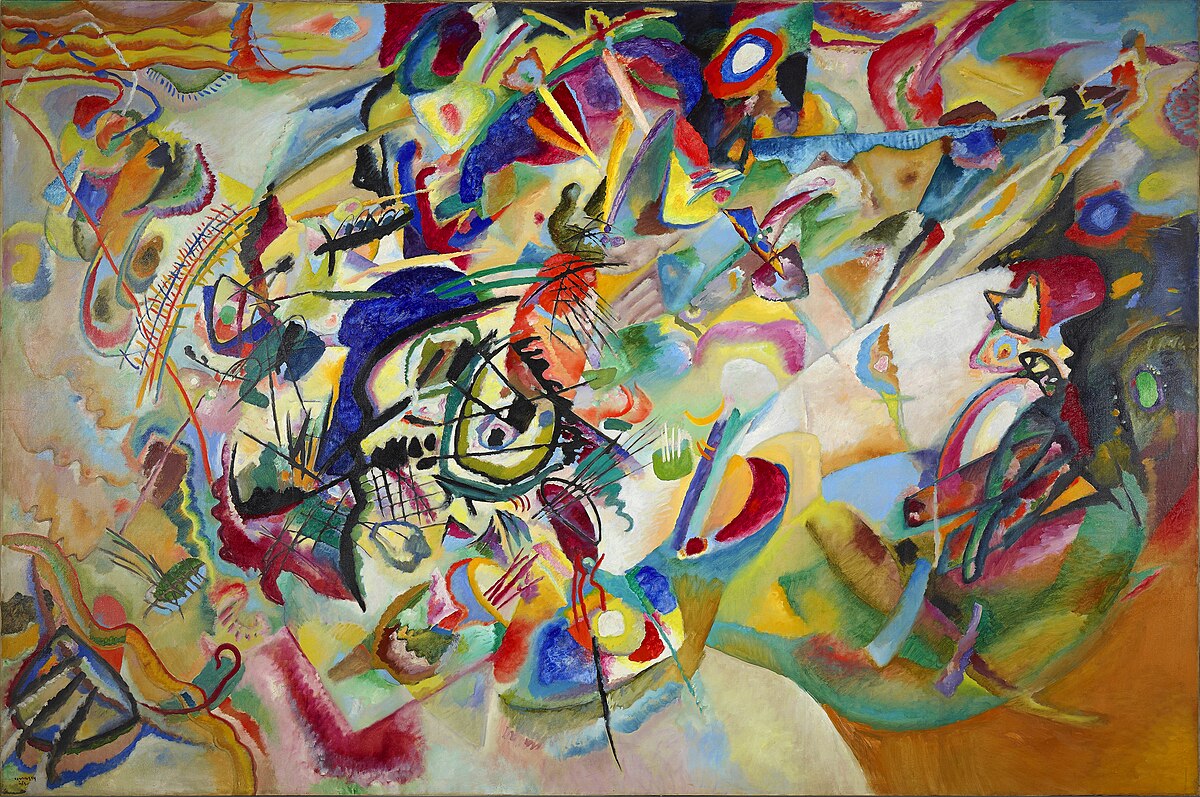}
    \caption{Figure 4.5: Wassily Kandinsky’s \emph{Composition VII} (1913).}
\end{subfigure}
\hfill
\begin{subfigure}[t]{0.43\textwidth}
    \centering
    \includegraphics[width=\textwidth]{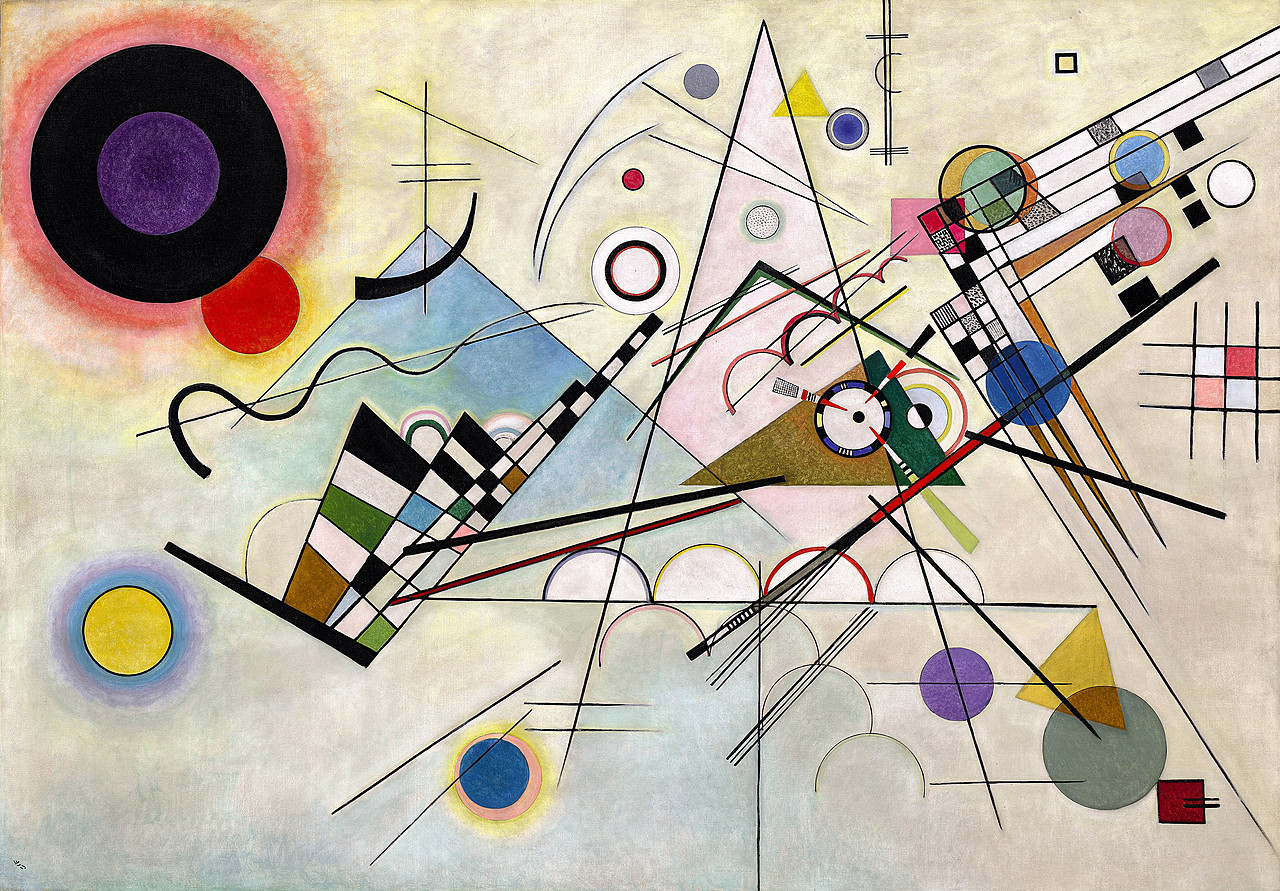}
    \caption{Figure 4.6: Wassily Kandinsky’s \emph{Composition VIII} (1923).}
\end{subfigure}

\caption{Comparison of Kandinsky’s \emph{Composition VII} and \emph{Composition VIII}.}
\end{figure}

\paragraph{Composition VII as a Nonlinear Dynamical Flow:}
The sweeping arcs, curves, and color gradients of \emph{Composition VII} can be modeled as the evolution 
of an initial geometric field $\phi_0 \in H$ under a nonlinear flow 
\[
\phi(t) = \Phi_t(\phi_0),
\]
where $\{\Phi_t\}$ is a semigroup of diffeomorphic transformations representing stretching, shearing, and 
curvilinear advection on $\Omega$.  
These operators are analogous to those arising in fluid dynamics or nonlinear PDEs, generating continuous 
deformations rather than discrete geometric orbits.  
The painting's dense layering corresponds to superpositions 
$\sum_{k=1}^{N} \phi(t_k)$ at different times, so that \emph{Composition VII} may be viewed as a 
time-integrated field of nonlinear operator actions, producing the visual analogue of a dynamical system 
with high-dimensional trajectories.

\paragraph{Composition VIII as an Affine Operator System.}
Kandinsky’s \emph{Composition VIII} (1923) may be interpreted as the realization of an operator system 
acting on a finite set of geometric generators.  
Let $\Omega \subset \mathbb{R}^{2}$ denote the canvas and consider the Hilbert space 
$H = L^{2}(\Omega;\mathbb{R}^{3})$, where the three components represent the RGB (or CMY) color channels.  
Fix a small collection of prototype shapes 
$\{\phi_j\}_{j=1}^{J} \subset H$, such as disks, line segments, triangles, and polygons, each given by the 
indicator of a region multiplied by a constant color vector.  
Let $\mathcal{G}$ be the semigroup generated by translations $T_a$, rotations $R_\theta$, and anisotropic 
scalings $S_\lambda$ on $\Omega$; each element $U \in \mathcal{G}$ acts unitarily (up to scaling) on $H$ via 
$(Uf)(x) = f(U^{-1}x)$.  

Then the visual content of \emph{Composition VIII} can be modeled as the finite superposition
\[
f(x) = \sum_{k=1}^{N} c_k \, U_k \phi_{j(k)}(x),
\]
where $c_k \in \mathbb{R}^{3}$ encodes color intensity and $U_k \in \mathcal{G}$ encodes the spatial 
placement, orientation, and size of each element.  
The orbit set $\{U_k \phi_{j(k)}\}_{k=1}^{N}$ forms a structured frame for $\operatorname{span}\{\phi_j\}$, 
so that geometric redundancy and overlap yield a stable representation of the composition.  
Kandinsky’s interplay of circles, radiating lines, and intersecting polygons corresponds, in this model, to 
the interaction of commuting and non-commuting operators in $\mathcal{G}$, producing visual analogues of 
spectral decomposition (concentric circles as ``eigenmodes'') and interference (crossing lines as superposed 
directions).  
Thus, \emph{Composition VIII} may be viewed as an explicit example in which an abstract operator system on 
a geometric basis in $H$ generates a complex, yet mathematically controlled, field of visual forms.
\paragraph{Comparative Mathematical Structure.}
Taken together, \emph{Composition VII} and \emph{Composition VIII} form complementary mathematical models: 
the former is governed by nonlinear, continuous deformation and high entropy in shape space, whereas the 
latter embodies discrete affine symmetries, group orbits, and frame-theoretic redundancy.  
Both works can therefore be situated within a unified operator-theoretic framework, highlighting how 
Kandinsky’s abstraction arises from the controlled action of transformation systems on geometric and 
color fields in $H$.

\subsubsection{Spectral, Algebraic, and Topological Structures in Kandinsky’s Compositions}

Beyond their visual differences, both \emph{Composition VII} and \emph{Composition VIII} can be interpreted 
through deeper spectral and algebraic structures.  
If $U$ denotes a geometric transformation in either the nonlinear flow $\{\Phi_t\}$ or the affine semigroup 
$\mathcal{G}$, then the induced operator on $H$, defined by $(Uf)(x)=f(U^{-1}x)$, generates a 
representation of a transformation group on the visual field.  
The circular motifs in \emph{Composition VIII} behave as approximate eigenfunctions of radial operators, 
while the sweeping curves in \emph{Composition VII} encode the integral curves of vector fields on $\Omega$, 
suggesting an underlying topological flow.  
Intersections of lines and arcs form nodes of operator interaction, analogous to the superposition of 
eigenmodes in harmonic analysis.  
Thus, both compositions reveal a spectral--topological organization in which Kandinsky’s abstract forms arise 
from the algebraic action of transformation systems, framing his art as an interplay between geometry, 
motion, and operator theory.

\section{The Fibonacci Sequence in Nature and Art}
In nature, the Fibonacci sequence appears in various forms, such as the arrangement of leaves on a stem, the branching of trees, and the pattern of seeds in a sunflower. These natural patterns have inspired artists to incorporate the Fibonacci sequence into their work, creating a bridge between the organic world and human creativity.\\

From an artistic standpoint, the Fibonacci sequence represents the mathematical essence of growth, rhythm, and proportion. The ratio $\varphi$ has historically guided architects and painters — from the Parthenon’s façade to Leonardo da Vinci’s \emph{Vitruvian Man}. In mathematical painting, using Fibonacci-based grids, spiral compositions, or scaling factors introduces a natural sense of order within abstraction. Each iteration visually mirrors the recursive logic of the sequence: expansion without distortion, harmony through self-reference.

Thus, the Fibonacci sequence is not only a numerical curiosity but a foundational bridge between mathematics, geometry, and art — a universal rhythm that defines both structural beauty and creative balance.

\section*{Mathematical Foundations of the Fibonacci Sequence}

The Fibonacci sequence, defined recursively by  
\[
x_{n+1} = x_n + x_{n-1}, \qquad x_1 = 0,\; x_2 = 1,
\]
gives rise to the well-known progression
\[
0,\,1,\;1,\;2,\;3,\;5,\;8,\;13,\ldots
\]
whose influence extends far beyond number theory. Fibonacci numbers appear naturally in phyllotaxis, biological growth patterns, the proportions of shells and flowers, architectural design, and numerous artistic compositions. The recursive nature of the sequence yields ratios converging to the Golden Ratio $\varphi = \frac{1+\sqrt{5}}{2}$, generating the characteristic scaling patterns and spirals that define many self-similar structures in nature.  

\begin{figure}[h!]
\centering
\includegraphics[width=0.45\textwidth]{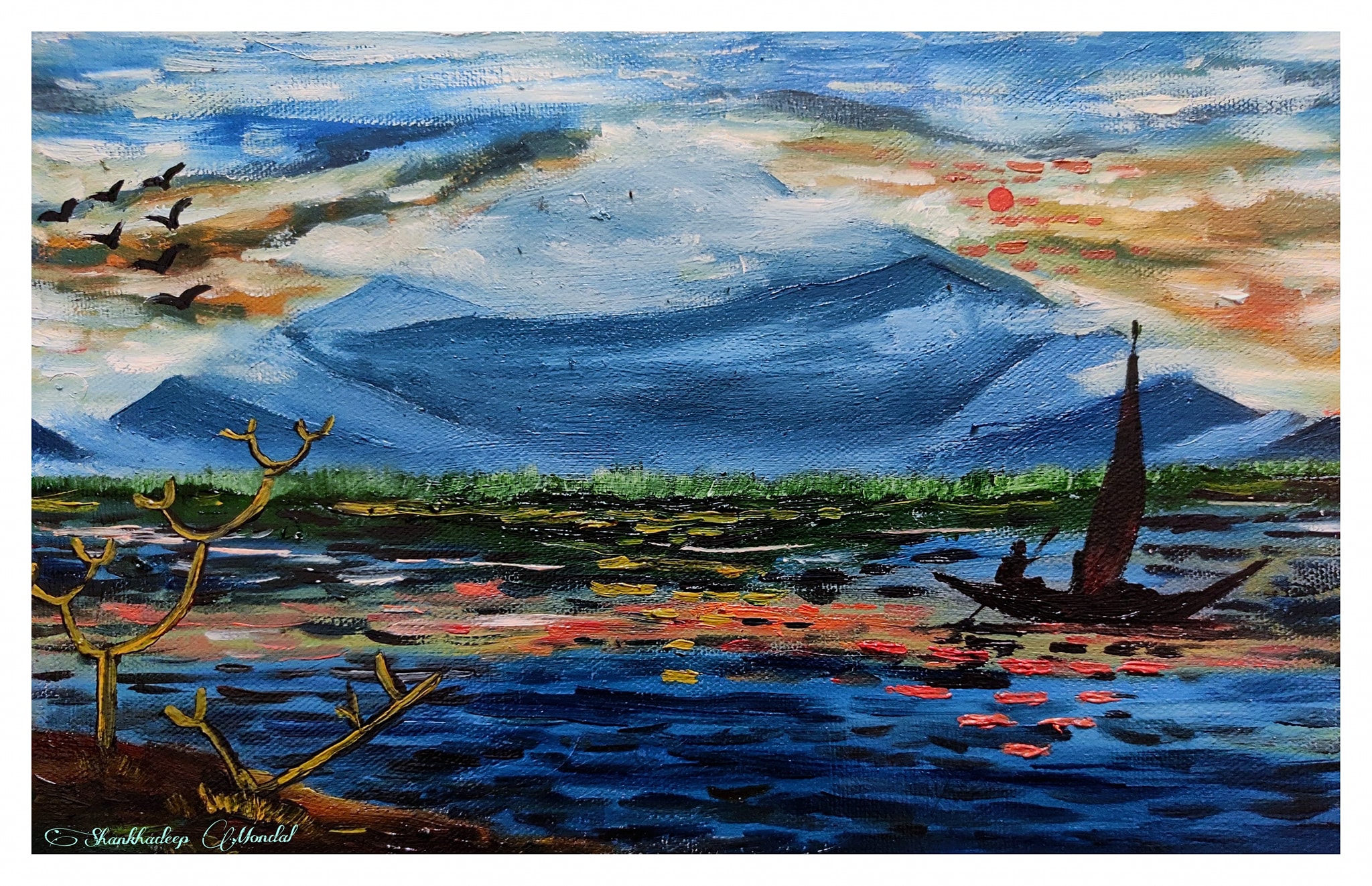}
\caption*{\textbf{Figure 5.1} \emph{Where Mathematics Meets Art} by Shankhadeep Mondal.}
\end{figure}

In the painting \emph{``Where Mathematics Meets Art''}, the Fibonacci sequence serves as the underlying generative structure for the entire composition. The clustering of birds, the reflections of sunlight on the water, and the distribution of clouds follow Fibonacci-like groupings, with visual elements increasing in quantities and separations according to $1,2,3,5,8$. The foreground trees are constructed using a Fibonacci branching scheme, in which each new branch emerges as the sum of the two preceding directions, mirroring the recurrence relation $x_{n+1}=x_n+x_{n-1}$.

The geometry of the mountains forms the central mathematical motif of the work. When the mountain peaks are interpreted as points on an $X,Y$ plane, the horizontal and vertical distances between them follow the proportional pattern  
\[
1 : 1 : 2 : 3 : 5 : 8,
\]
revealing a landscape structured by the Fibonacci sequence. Connecting the endpoints of these mountain divisions yields a total count that corresponds to the next Fibonacci number, visually expressing the recursive nature of the sequence within the terrain. Similar proportional logic appears in the contours of the landmasses and in the placement of the boat, whose orientation lies along a Fibonacci-guided visual trajectory.  

Together, these elements demonstrate how the Fibonacci sequence can govern an artwork not only symbolically but structurally, shaping spatial rhythm, compositional balance, and the emergent harmony of the scene. The painting thus illustrates how mathematical recursion can seamlessly integrate with natural imagery, forming a bridge between numerical abstraction and organic representation.

\begin{figure}[h!]
\centering
\includegraphics[width=0.45\textwidth]{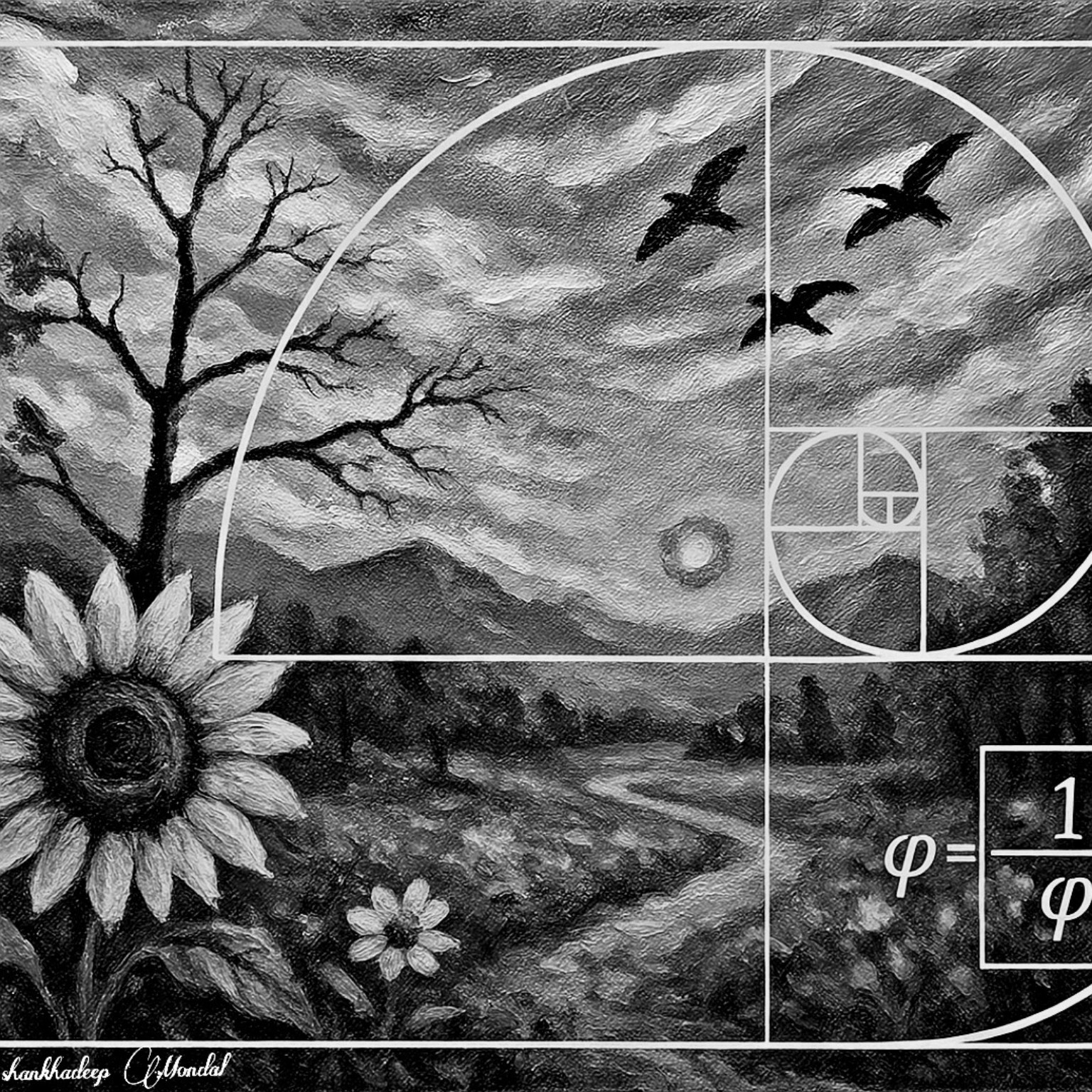}
\caption*{\textbf{Figure 5.2} \emph{Symmetry of the Infinite} by Shankhadeep Mondal.}
\end{figure}

\noindent
The second painting, \emph{``Symmetry of the Infinite''}, presents a monochromatic landscape in which the Fibonacci sequence and the Golden Ratio $\varphi$ are explicitly embedded into the geometric structure of the composition. The sunflower exemplifies Fibonacci phyllotaxis, where the number of spirals in each direction reflects consecutive Fibonacci numbers, while the branching of the tree follows the recurrence $x_{n+1}=x_n+x_{n-1}$. The arrangement of clouds, the flight of the birds, and the layering of the mountains align with the expanding Golden Spiral generated by nested Golden Rectangles. Even the curvature of the path and the placement of the sun correspond to $\varphi$-scaled divisions across the canvas.

\subsubsection*{Solitude in the Geometry of Dreams}

\begin{figure}[h!]
\centering
\includegraphics[width=0.50\textwidth]{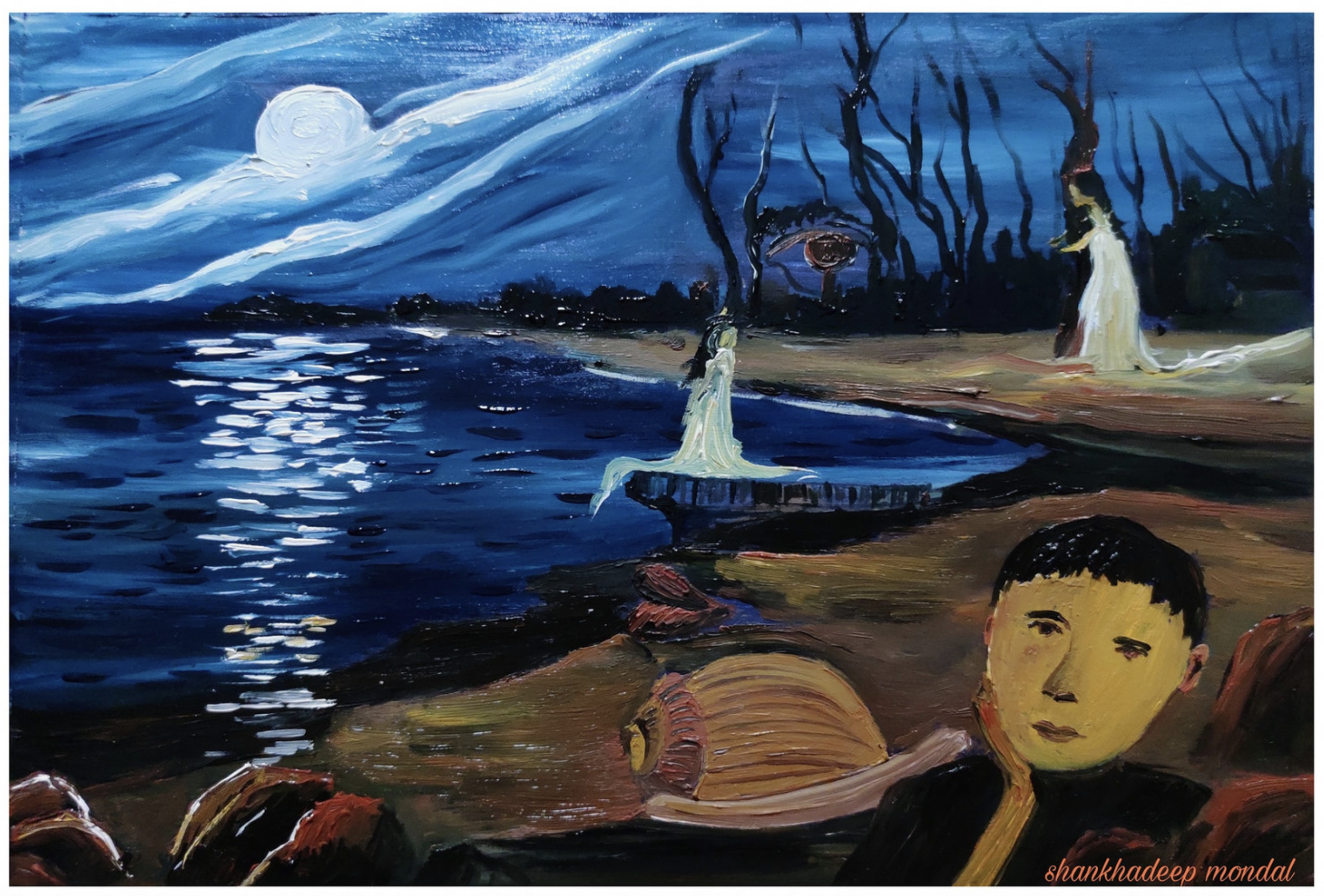}
\caption*{\textbf{Figure 5.3} \emph{Solitude in the Geometry of Dreams} by Shankhadeep Mondal.}
\end{figure}

\noindent
This surrealist painting blends emotional solitude with a subtle mathematical undercurrent. The moonlit sky and its reflection on the water form a pair of approximately symmetric curves, evoking the geometric idea of reflection across a horizontal axis. The elongated cloud bands resemble smoothed waveforms, suggestive of oscillatory behaviour that can be modeled by continuous functions on $\mathbb{R}$. The placement of the figures and objects along the shoreline follows a loose recursive spacing, echoing the structure of a non-uniform sequence $\{x_n\}$ whose increments gradually diminish, creating a visual rhythm of contraction toward the distant horizon. Through these elements, the painting employs mathematical suggestions—symmetry, oscillation, recursion—to deepen its atmosphere of introspective mystery under the infinite night.

\section{Cubism, Higher-Dimensional Mathematics, and Fractal Geometry}

Historically, the early cubists were exposed to higher-dimensional geometry through the influential work of Esprit Jouffret, whose 1903 treatise on four-dimensional polytopes provided explicit diagrammatic projections of hypercubes, hypersimplices, and polytopal unfoldings. Jouffret’s diagrams, in turn, were grounded in the differential geometric developments of Gauss, Riemann, and Poincaré, who challenged the axioms of Euclidean space and introduced curvature, higher-dimensional manifolds, and non-Euclidean structures. These ideas, communicated to Picasso, Braque, and their contemporaries by Maurice Princet, the so-called "mathematician of Cubism," precipitated the abandonment of the single central projection in favor of a multi-projectional, multi-dimensional mapping.

\subsection{Cubism, Higher Mathematics, and Fractal Geometric Structure}

The evolution of Cubism in the early twentieth century coincided with major developments in geometry and analysis, particularly the work of Riemann, Poincar\'e, and Jouffret on non-Euclidean manifolds, higher-dimensional embeddings, and projection theory. Through Maurice Princet, Picasso and Braque were introduced to Jouffret's diagrams of four-dimensional polytopes, which illustrated the planar shadows of hypercubes and other multi-dimensional bodies. In mathematical terms, a cubist representation may be modeled as a multi-projection operator system
\[
\mathcal{C}(X)=\sum_{\theta\in\Theta} w(\theta)\, P_{\theta}(X),
\]
where each \(P_{\theta}:\mathbb{R}^{3}\to\mathbb{R}^{2}\) denotes a distinct geometric viewpoint. This replaces the classical single-point linear projection with a superposition of multiple incompatible orientations, effectively simulating the projection of an embedded object \(X\hookrightarrow \mathbb{R}^{4}\) onto the plane.

The decomposition of visual objects into planar facets parallels a frame expansion in the Hilbert space \(H=L^{2}(\Omega)\), where the image \(f\) is represented as
\[
f=\sum_{k}\langle f,g_{k}\rangle\phi_{k},
\]
with \(\{\phi_{k}\}\) serving as geometric atoms such as triangles, quadrilaterals, or curved patches---forming a redundant and stable frame. Such redundancy anticipates robustness under erasure, a key feature in modern harmonic analysis and frame theory.

Moreover, many cubist surfaces exhibit hierarchical, self-affine subdivision analogous to fractal Iterated Function Systems (IFS). For affine contractions \(T_{i}(x)=A_{i}x+b_{i}\), a fractal set \(F\) satisfies
\[
F=\bigcup_{i=1}^{N} T_{i}(F),
\]
and its Hausdorff dimension \(d\) is determined by Moran's equation \(\sum_{i=1}^{N} r_{i}^{\,d}=1\), where \(r_{i}\) are the scaling ratios. Cubist compositions often emulate this self-affine structure, producing multi-scale geometric repetition and visual complexity. 

Thus, Cubism may be rigorously interpreted as the intersection of higher-dimensional projection geometry, fractal self-similarity, and operator-theoretic decomposition, embedding artistic abstraction within a unified and deeply mathematical framework.

\subsubsection{Cubism and Frame-Theoretic Decomposition}

The decomposition of objects into overlapping geometric facets in Cubism parallels a frame expansion in 
the Hilbert space $H=L^{2}(\Omega)$.  
If $\{\phi_{k}\}$ denotes the collection of planar fragments appearing in a cubist composition, the image 
may be represented as
\[
f = \sum_{k} \langle f, g_{k} \rangle \phi_{k},
\]
where $\{g_{k}\}$ is a dual frame.  
The inherent redundancy of frames mirrors the cubist principle of multiple simultaneous viewpoints, 
yielding stable reconstructions even under partial erasure of visual information.

\subsubsection{Cubism as a Multi-Projection Operator System}

A central mathematical interpretation of Cubism arises from replacing the classical single-view projection 
\(P:\mathbb{R}^{3}\to\mathbb{R}^{2}\) with a superposition of multiple projections.  
Given a family of viewpoint-dependent operators 
\(\{P_{\theta}:\mathbb{R}^{3}\to\mathbb{R}^{2}\}_{\theta\in\Theta}\) and weights \(w(\theta)\ge0\) with 
\(\sum_{\theta\in\Theta} w(\theta)=1\), a cubist representation of an object \(X\subset\mathbb{R}^{3}\) is modeled by
\[
\mathcal{C}(X)=\sum_{\theta\in\Theta} w(\theta)\, P_{\theta}(X).
\]
This multi-projection operator captures the simultaneous appearance of incompatible viewpoints, 
effectively simulating a higher-dimensional embedding followed by a composite planar projection.

\subsubsection{Cubism and Fractal Geometry: Self-Similarity and Iterated Operators}

The recursive subdivision of planes in many cubist works parallels the action of self-similar and 
self-affine Iterated Function Systems (IFS).  
Let $\{T_i\}_{i=1}^{N}$ be affine maps on $\mathbb{R}^2$, $T_i(x)=A_i x + b_i$, where each $A_i$ is contractive.  
The associated fractal set $F$ is the fixed point of Hutchinson's operator,
\[
F = \bigcup_{i=1}^{N} T_i(F),
\]
mirroring the repeated faceting and geometric refinement in analytical cubism.  
Such structures exhibit non-integer dimensionality, with Hausdorff dimension $d$ determined by 
Moran's equation $\sum_{i=1}^{N} r_i^{\, d}=1$, where $r_i$ are the contraction ratios.  
Thus, cubist decomposition can be viewed as a visual analogue of fractal self-similarity and 
iterated operator dynamics.
\subsubsection{Cubism and Fractal Measures}

The geometric decomposition in cubist compositions can be modeled using self-affine fractal measures.  
Let $\mu$ be a probability measure satisfying the invariance relation
\[
\mu = \sum_{i=1}^{N} p_i\, \mu \circ T_i^{-1},
\]
where $\{T_i\}$ are affine contractions and $\{p_i\}$ are weights.  
Such measures appear in the analysis of fractal frames and wavelets on self-similar sets, and the 
distribution of planar facets in cubist works resembles the support of $\mu$.  
This correspondence links the multiscale visual density of cubism with the harmonic and geometric 
properties of fractal measures.

\subsection{Cubism as a Higher-Dimensional Embedding Problem}

Cubist simultaneity may be interpreted mathematically as the projection of a higher-dimensional object 
onto the plane.  
Let $E:X\hookrightarrow \mathbb{R}^{n}$, $n\ge 4$, be an embedding of a visual object, and let 
$P_{\theta}:\mathbb{R}^{n}\to\mathbb{R}^{2}$ denote viewpoint-dependent projections.  
A cubist image is then modeled by the composite map
\[
\mathcal{C}(X)=\sum_{\theta\in\Theta} w(\theta)\, P_{\theta}(E(X)),
\]
where $w(\theta)$ are nonnegative weights.  
This framework interprets cubist geometry as the superposition of multiple $n$-dimensional shadows, 
reflecting the collapse of higher-dimensional structure into a planar representation.

\subsection{Example: A Cubist--Fractal Geometric Composition}
\begin{figure}[h!]
\centering
\includegraphics[width=0.45\textwidth]{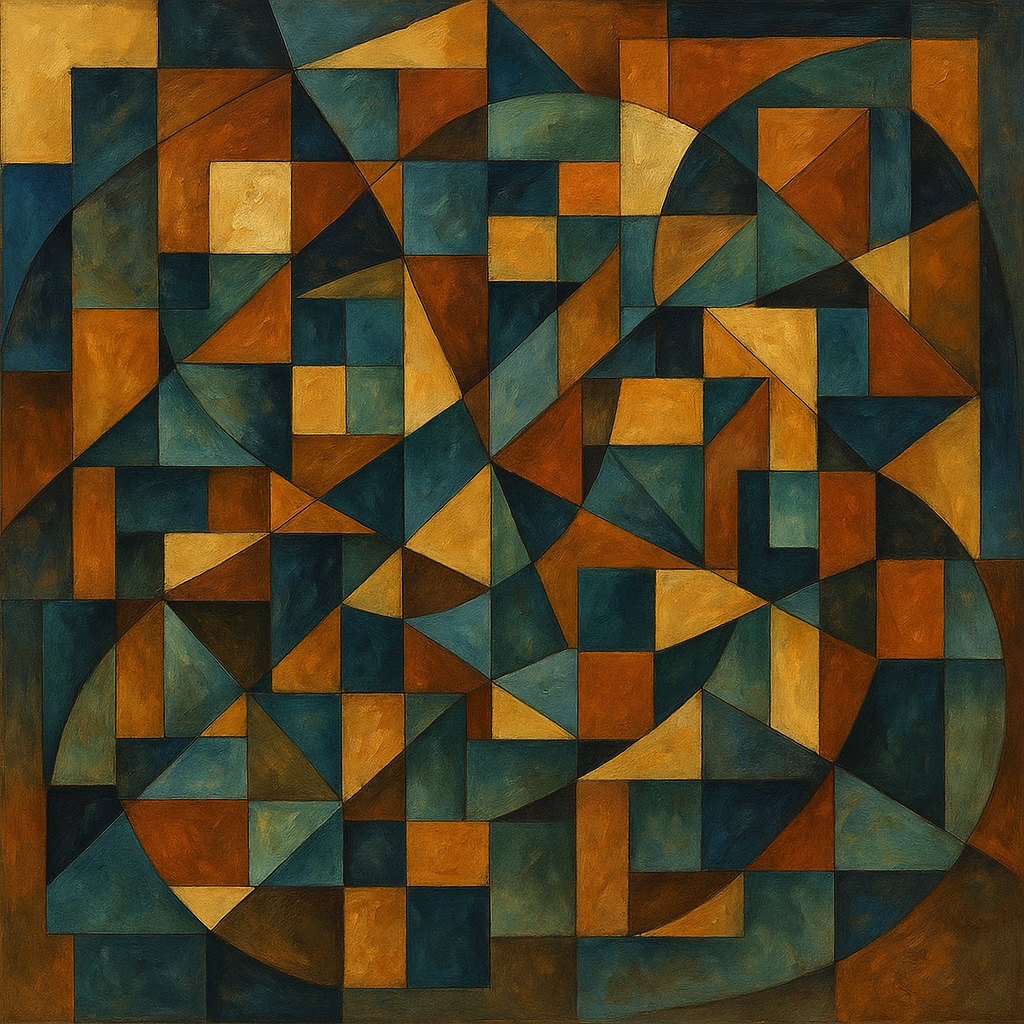}
\caption*{\textbf{Figure 6.1} \emph{Symphony of Shapes} by Shankhadeep Mondal.}
\end{figure}
The painting in Fig.6.1 illustrates a cubist composition constructed from  overlapping triangles, quadrilaterals, and curved segments, each representing a distinct projection of an underlying three-dimensional object.  The visual field may be modeled as a function 
$f \in L^{2}(\Omega)$ decomposed into geometric atoms $\{\phi_{k}\}$, yielding the frame expansion 
$f=\sum_{k}\langle f,g_{k}\rangle\phi_{k}$.  
The recursive subdivision of the canvas into smaller self-affine regions mirrors an Iterated Function 
System $\{T_i\}$ with invariant set $F=\bigcup_i T_i(F)$, producing multi-scale structure reminiscent of 
fractal measures.  
The curved arcs correspond to projections of a higher-dimensional embedding 
$X\hookrightarrow\mathbb{R}^{4}$, while the planar facets represent the shadows $P_{\theta}(X)$ arising 
from multiple viewpoints.  
Thus, the image serves as a concrete example of how cubist geometry integrates multi-projection 
operators, fractal self-similarity, and frame-theoretic redundancy within a single mathematical 
representation.

\begin{figure}[h!]
\centering
\includegraphics[width=0.50\textwidth]{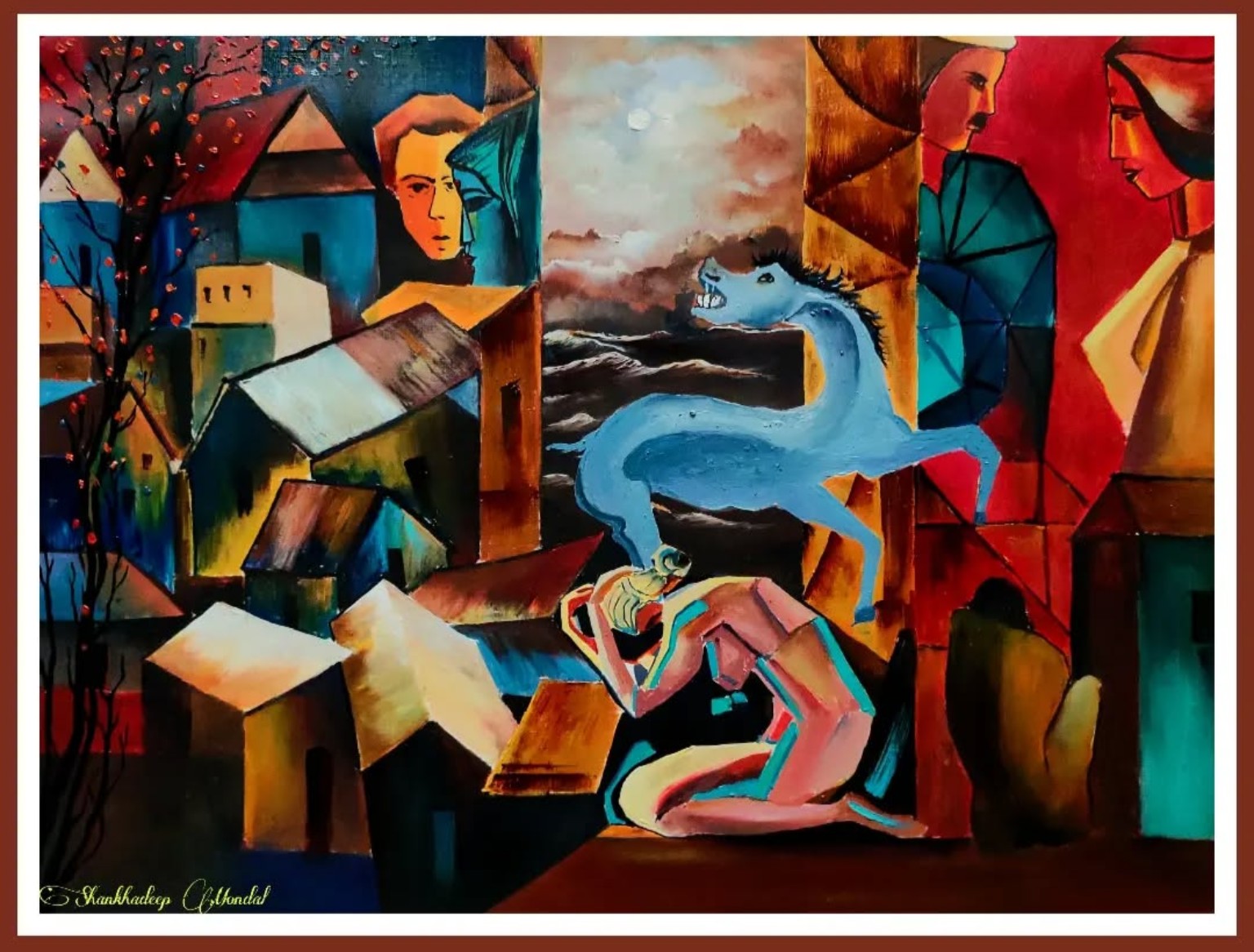}
\caption*{\textbf{Figure 6.2} \emph{Life and War} -- A Cubist Composition by Shankhadeep Mondal.}
\end{figure}

\noindent
\emph{Life and War} employs cubist geometry to encode emotional tension through mathematical fragmentation. The composition is built from intersecting polygons, affine transformations, and planar subdivisions, reflecting the breakdown of form under conflict. The figures and architectural blocks are decomposed into triangular and quadrilateral facets, each corresponding to local linear transformations of an underlying shape space. The distortions of the horse and the kneeling figure arise from non-uniform scalings and rotations, while the surrounding faces emerge through symmetry-breaking across adjacent planes. This geometric disassembly---a hallmark of cubism---creates a visual topology in which multiple perspectives coexist simultaneously, revealing war as a nonlinear disruption of the ordinary metric structure of life.

\section{Techniques for Incorporating the Fibonacci Sequence in Painting}

Artists often incorporate the Fibonacci sequence through geometric tools such as the golden rectangle and the golden spiral. A common technique is to construct a Fibonacci tiling---a sequence of adjacent squares whose side lengths follow $1,1,2,3,5,8,\ldots$---and then draw circular arcs across opposite corners to create the golden spiral. This spiral is frequently used to structure the composition, establish focal points, and guide the viewer’s eye naturally across the artwork.

\subsection*{Case Studies}
Leonardo da Vinci’s \emph{Vitruvian Man} demonstrates golden-ratio proportions in the human body;  
Dalí’s \emph{The Sacrament of the Last Supper} uses golden rectangles to shape the spatial geometry;  
and many contemporary artists employ Fibonacci grids and digital rendering to compose mathematically balanced designs.

\subsection*{Algorithm for Constructing a Fibonacci-Based Composition}

\begin{algorithm}[h!]
\caption{Generating a Fibonacci-Guided Painting Layout}
\begin{algorithmic}[1]
\State Choose an initial square of side length $1$.
\State Append another square of side length $1$ adjacent to the first.
\For{$n = 3$ to $N$}
    \State Construct a new square of side length $x_{n} = x_{n-1} + x_{n-2}$,
           attaching it to the existing rectangle in a counterclockwise pattern.
\EndFor
\State Draw quarter-circle arcs inside each square to form the golden spiral.
\State Place major focal elements (figures, objects, light sources) along points of high curvature on the spiral.
\State Align secondary structures (mountains, trees, buildings) with the edges of the Fibonacci rectangles.
\State Adjust color gradients and shadow lines to follow the spiral’s direction to enhance visual flow.
\State Render final details while preserving the global structure defined by the Fibonacci grid.
\end{algorithmic}
\end{algorithm}

\section{Convergence of Cubist Structure and Fibonacci Proportion}

The relationship between Cubism, the Fibonacci sequence, and the Golden Ratio reveals a compelling fusion of geometric abstraction and mathematical proportion. Although early Cubist painters did not explicitly claim mathematical intentions, many works from this period exhibit spatial harmonies aligned with $\varphi=\frac{1+\sqrt{5}}{2}$ and with recursive geometric division. The fractured planes, shifting perspectives, and affine transformations characteristic of Cubism naturally lend themselves to proportional systems derived from golden rectangles and Fibonacci tilings.

A striking historical example of this convergence is Marcel Duchamp’s \emph{Nude Descending a Staircase (No.~2)} (1912). Its sequential fragmentation creates a rhythmic progression reminiscent of the Fibonacci recurrence $x_{n+1}=x_n+x_{n-1}$, while several bounding regions align closely with golden-ratio subdivisions. Similarly, in Picasso’s \emph{Girl with a Mandolin} (1910) and Braque’s \emph{Violin and Candlestick} (1910), the underlying arrangement of planes corresponds to proportional grids that echo golden-ratio divisions, demonstrating how Cubist compositions often embodied mathematical balance even when not explicitly acknowledged.

Modern and contemporary artists have extended this dialogue. Mario Merz, for example, integrated explicit Fibonacci numbers and spirals into structurally fragmented works, showing how recursive growth patterns can coexist with Cubist geometry. Digital and algorithmic art further expands these possibilities: computational tools now generate Cubist compositions governed by Fibonacci spirals, golden grids, and recursive subdivision algorithms, revealing that the connection between Cubism and mathematical proportion is both historically rooted and creatively evolving.

Together, these examples show that Fibonacci structure and golden-ratio geometry provide a powerful lens for understanding Cubist form. Mathematical growth laws and proportional harmonies can integrate seamlessly with the visual logic of fragmentation, producing artworks that operate simultaneously on analytical, aesthetic, and conceptual levels.

\section{Discussion}

The use of the Fibonacci sequence and the golden ratio in painting illustrates a compelling intersection between mathematical structure and artistic creativity. Historically, artists have employed Fibonacci-based proportions---such as golden rectangles and spirals---to achieve visual balance, rhythmic spacing, and coherent compositional flow. These constructions often align with naturally occurring patterns, reinforcing a sense of harmony that viewers may intuitively perceive. Although psychological studies debate whether the golden ratio holds universal aesthetic appeal, it remains clear that such mathematical frameworks can subtly influence the viewer’s perception of order and visual comfort.

In contemporary art, digital tools and algorithmic methods have broadened the possibilities for incorporating Fibonacci structures. Through computational geometry, artists can now generate complex recursive patterns and precise proportional grids that would be difficult to construct manually. These modern techniques highlight how the Fibonacci sequence continues to inspire innovative approaches to visual design. Ultimately, while mathematical structure contributes to compositional clarity and aesthetic resonance, its impact varies with individual perception, emphasizing that the Fibonacci sequence functions as a flexible guiding principle rather than a rigid artistic rule.

\section{Conclusion and Future Perspectives}

The Fibonacci sequence and the Golden Ratio continue to exert a profound influence on the theory and practice of painting, shaping visual structure from Renaissance proportion studies to contemporary algorithmic art. Their recursive and proportional properties provide artists with a mathematical foundation for balance, harmony, and visual flow, while also bridging natural growth patterns and artistic intuition. As shown throughout this work, the Fibonacci sequence informs not only geometric layouts and compositional grids but also symbolic interpretation, aesthetic resonance, and psychological engagement.

Looking ahead, several promising research directions emerge. One avenue concerns empirical studies of how viewers perceive Fibonacci-based compositions: psychophysical experiments, eye-tracking, and neuroaesthetic methods may clarify whether golden-ratio structures genuinely guide attention or evoke distinct emotional responses. Another direction lies in computational approaches, including recursive subdivision algorithms, parametric design tools, and machine-learning models for generating Fibonacci-guided artworks. The integration of Fibonacci geometry into virtual and augmented reality environments further opens possibilities for immersive, three-dimensional visual experiences shaped by recursive scaling.

Interdisciplinary work may also explore connections between Fibonacci geometry and other artistic systems such as Cubism, fractal art, or dynamical-pattern generation. In this context, Sir Roger Penrose’s aperiodic tilings provide a particularly rich framework. Constructed from non-periodic arrangements governed by golden-ratio proportions, Penrose tilings offer a mathematically precise yet visually poetic means of creating compositions that exhibit global order without translational repetition, suggesting fertile ground for hybrid visual systems.

Beyond Fibonacci structures, the Lucas sequence defined by the same recurrence with different initial conditions offers a closely related but distinct geometric framework whose applications in painting remain largely unexplored. Although Lucas numbers also converge asymptotically to the golden ratio, their differing combinatorial and scaling properties suggest alternative proportional systems, tilings, and growth patterns. At present, no systematic theory or established artistic practice incorporates Lucas-sequence geometry into visual composition, making this an open and promising direction for future research.

In summary, recursive sequences provide a powerful mathematical grammar for artistic composition. The Fibonacci sequence remains not only a historical emblem of mathematical beauty but also a dynamic source of contemporary artistic inquiry. Its continued exploration alongside related sequences such as the Lucas sequence promises to expand both theoretical understanding and creative possibilities in mathematically informed visual art.

\section{Acknowledgment}
\noindent The authors are grateful to the Mohapatra Family Foundation and the College of Graduate Studies of the University of Central Florida for their support during this research.

\section*{References}

\begin{itemize}

    \bibitem{livio2002golden}
    Livio, M. (2002). \emph{The Golden Ratio: The Story of Phi, the World's Most Astonishing Number}, Broadway Books.

    \bibitem{} Dunlap, R. A. (1997). \emph{The Golden Ratio and Fibonacci Numbers}, World Scientific.

    \bibitem{huntley1972divine} Huntley, H. E. (1972). \emph{The Divine Proportion: A Study in Mathematical Beauty}, Dover Publications.

    \bibitem{ghyka1977geometry}
    Ghyka, M. (1977). \emph{The Geometry of Art and Life}, Dover Publications.

   \bibitem{jensen2002mathematics}
   H J Jensen, (2002). \emph{Mathematics and painting}, Interdisciplinary science reviews.
   
   \bibitem{} P. R. Halmos,  \emph{MATHEMATICS AS A CREATIVE ART}, The Scientific Research Honor Society, Winter 1968, Vol. 56, No. 4 (Winter 1968), pp. 375-389
   
    \bibitem{doczi1994power}
    Doczi, G. (1994). \emph{The Power of Limits: Proportional Harmonies in Nature, Art, and Architecture}. Shambhala.

    \bibitem{field1997invention}
    Field, J. V. (1997). \emph{The Invention of Infinity: Mathematics and Art in the Renaissance}, Oxford University Press.

    \bibitem{henderson1983fourth}
    Henderson, L. D. (1983). \emph{The Fourth Dimension and Non-Euclidean Geometry in Modern Art}, Princeton University Press.

    \bibitem{} Lenoir, T. (1988). ``Aesthetics of Complexity: Duchamp, Futurism, and the Fourth Dimension,'' \emph{Configurations}, 1(2), 201–243.

    \bibitem{} Bradley, J. J.,  Bradley, J. (2015). ``Mathematics and Cubism: Proportions and Structures in the Works of Picasso and Braque,'' \emph{Journal of Mathematics and the Arts}, 9(3), 113–129.

    \bibitem{kaplan2013algorithmic}
    Kaplan, C. S. (2013). ``Algorithmic Beauty in Mathematical Art,'' \emph{Leonardo}, 46(1), 12–22.

    \bibitem{} Penrose, R. (1974). ``The Role of Aesthetics in Pure and Applied Mathematical Research,'' \emph{Bulletin of the Institute of Mathematics and its Applications}, 10, 266–271.

    \bibitem{} Penrose, R. (1974). ``The Riddle of the Tiling,'' \emph{Scientific American}, 236(1), 112–123.

    \bibitem{} Merz, M. (1980). \emph{Fibonacci and the Spiral in Contemporary Art}, Arte Povera Publications.

    \bibitem{} Zeki, S. (1999). \emph{Inner Vision: An Exploration of Art and the Brain}, Oxford University Press.

    \bibitem{} Zeki, S. (2013). ``Clues to the Biological Bases of Mathematical Beauty,'' \emph{Frontiers in Human Neuroscience}, 7, 1–12.

    \bibitem{} Weyl, H. (1952). \emph{Symmetry}, Princeton University Press.

    \bibitem{} Stewart, I. (2007). \emph{Why Beauty Is Truth: The History of Symmetry}, Basic Books.

    \bibitem{} Arnheim, R. (1974). \emph{Art and Visual Perception: A Psychology of the Creative Eye}, University of California Press.

    \bibitem{} Ostwald, M. J.,  Vaughan, J. (2016). ``Fibonacci Numbers, Proportional Systems, and Architectural Composition,'' \emph{Nexus Network Journal}, 18(2), 225–247.
    
\end{itemize}

\section{Biographical Information}

\textbf{Dr. Shankhadeep Mondal} received his M.S. and Ph.D. in Mathematics from IISER Thiruvananthapuram. His research focuses on harmonic analysis, erasure-robust frames, and spectral optimization. He also explores the intersection of mathematics and visual art, with over 2500 paintings. He is currently a postdoctoral researcher at the University of Central Florida.\\

\textbf{Prof. Ram Narayan Mohapatra}  is a mathematician at the University of Central Florida, with previous appointments in Canada and Lebanon. He has authored over 250 research papers and two books, and edited eight conference proceedings. His current research interests include approximation theory, equilibrium problems, tensors, and frames in Hilbert $C^*$modules.

\end{document}